\documentclass[12pt]{article} 
\usepackage{amsmath, amssymb} 
\usepackage{amsthm} 
\usepackage{enumerate} 

\newtheoremstyle{plainsl}%
	{\topsep}
	{\topsep}
	{\slshape} 
	{}
	{\normalfont\bfseries}
	{.}
	{ }
	{}

\swapnumbers

{\theoremstyle{plainsl}
\newtheorem{theorem}{Theorem}[section]
\newtheorem{lemma}[theorem]{Lemma}
\newtheorem{corollary}[theorem]{Corollary}}
{\theoremstyle{remark}
}
\newtheorem{problem}[theorem]{Problem}

\newcommand\cref[1]{Corollary~\ref{cor:#1}}

\renewcommand\proof{\noindent\textsl{Proof. }}
\newcommand\sqr[2]{{\vbox{\hrule height.#2pt
    \hbox{\vrule width.#2pt height#1pt \kern#1pt
        \vrule width.#2pt}\hrule height.#2pt}}}
\renewcommand\qed{%
	\ifmmode\eqno\sqr53
	\else\nolinebreak\ \hfill\sqr53\medbreak\fi}

\numberwithin{equation}{section}

\newcommand\bfe[1]{{\textbf{e}}(#1)}
\newcommand\bs{\backslash}
\newcommand\copp{C_{\mathrm{opp}}}
\newcommand\ovy{{\comp Y}}
\newcommand\sqod[1]{\seq{#1}1d}
\newcommand\sqzr[1]{\seq{#1}0r}

\newcommand\eg[1]{\emph{Example~$#1$}\quad}

\newcommand\al{\alpha}
\newcommand\be{\beta}
\newcommand\de{\delta}

\newcommand\ga{\gamma}
\newcommand\Ga{\Gamma}

\newcommand\cD{{\mathcal D}}

\newcommand\comp[1]{{\mkern2mu\overline{\mkern-2mu#1}}}
\newcommand\diff{\mathbin{\mkern-1.5mu\setminus\mkern-1.5mu}}
\newcommand\sbs{\subseteq}

\newcommand\seq[3]{#1_{#2},\ldots,#1_{#3}}

\newcommand\aut[1]{{\rm Aut}(#1)}

\DeclareMathOperator{\Sym}{Sym}
\newcommand\sym[1]{\Sym(#1)}

\title{Completely Transitive Designs}

\author{ Chris D. Godsil and Cheryl E. Praeger}
\date{Manuscript from 1997; tidied a bit}

\begin{document}
\maketitle

\begin{abstract}
We view a design $\cD$ as a set of $k$-subsets of a fixed set $X$ of
$v$ points.  A $k$-subset of $X$ is at distance $i$ from $\cD$ if it
intersects some $k$-set in $\cD$ in $k-i$ points, and no subset in
more than $k-i$ points.  Thus $\cD$ determines a partition by distance
of the $k$-subsets of $X$.  We say $\cD$ is completely transitive if
the cells of this partition are the orbits of the automorphism group
of $\cD$ in its induced action on the $k$-subsets of $X$.  This paper
initiates a study of completely transitive
designs $\cD$. A classification is given of all examples for which 
the automorphism group is not primitive on $X$. In the primitive case 
the focus is on examples with the property that any two distinct $k$-subsets in
$\cD$ have at most $k-3$ points in common. Here a reduction is given 
to the case where the automorphism group is 2-transitive on $X$. 
New constructions are given by classifying all examples for some 
famlies of 2-transitive groups, leaving several unresolved cases.
\end{abstract}

\section{Introduction}

A partition $\pi=\sqod{C}$ of the vertex set of a graph $Y$ is {\sl
equitable} if, for each $i$ and $j$, the number of neighbours in
$C_j$ of a vertex in $C_i$ is determined by $i$ and $j$.  A typical
example is provided by the orbits of any group of automorphisms of
$Y$.  If $C\sbs Y$ then the \textsl{distance partition} with respect to
$C$ is the partition whose $i$-th cell consists of the vertices in $Y$
at distance $i$ from $C$, that is vertices whose minimum distance from
some vertex of $C$ is $i$. The maximum distance of a vertex in $Y$
from $C$ is the \textsl{covering radius} of $C$, and will usually be
denoted by $r$.  The distance partition relative to a subset $C$ is
not usually equitable, but when it is we call $C$ a \textsl{completely
regular} subset of $Y$. Moreover a completely regular subset $C$ is
said to be {\it completely transitive} if the cells of the distance
partition with respect to $C$ are the orbits of some group of
automorphisms of $Y$. Completely regular subsets are important in
coding theory and, more generally, in the theory of distance regular
graphs.  In particular, a regular graph is distance-regular if and
only if each vertex in it is a completely regular subset.  (For the
last claim see the proof of Theorem~2.2 in \cite{cg-sh}, for general
information on distance-regular graphs see \cite{bcn,g-bk} and for
more on equitable partitions \cite{g-bk}.)  Distance-transitive graphs
form an important and interesting subclass of distance-regular graphs.
In the present context we may may view these as vertex transitive
graphs with the property that each vertex is a completely transitive
subset.

The Johnson graph $J(v,k)$ is defined as follows.  Its vertices are
the $k$-subsets of a fixed subset of $v$ points and two $k$-subsets
are adjacent if they have exactly $k-1$ points in common.  Since
$J(v,k)$ and $J(v,v-k)$ are isomorphic, we will normally assume that
$2k\le v$.  The aim of this paper is to study completely transitive
subsets of $J(v,k)$.

Let $C$ be a subset of the vertex set $V(J(v,k))$ of $J(v,k)$.  We
have already defined the covering radius of $C$, but there are two
further parameters we need.  The \textsl{minimum distance} $\de$ of $C$
is simply the minimum distance between two vertices of $C$.  The \textsl{strength} 
of $C$ is the largest integer $t$ such that every
$t$-subset of the underlying $v$-set $X$ lies in the same number of
vertices of $C$. Thus, if $t\ge 1$, then $C$ is a $t$-design on $X$,
where here we are using $t$-design with its usual meaning, that is a
collection $C$ of $k$-subsets of $X$ with the property that each
$t$-subset of $X$ lies in the same number of elements of $C$. One
result of this paper (??) is that if $G$ is the automorphism group of
a completely transitive subset of $J(v,k)$ with $\de\ge3$ then $G$
must act 2-transitively on $X$, and hence $C$ has strength at least
two; that is $C$ is a 2-design.

Completely regular subsets of $J(v,k)$ are discussed at length in
\cite{mart,mar2,mar3}.

\section{Examples}

We shall describe the important known classes of completely transitive
designs in this section.  However we begin by making some general
observations about completely regular subsets.  First we note an
unpublished observation due to A. Neumaier.

\begin{lemma}
	\label{lem:CR} 
Let $C=C_0$ be a subset of a distance-regular
graph with distance partition $\sqzr{C}$.  Then $C_0$ is completely
regular if and only if $C_r$ is.
\end{lemma}

\proof
A simple induction argument shows that since the partition
$C_0,\ldots,C_r$ is equitable, any vertex in $C_i$ is joined by a path
of length $r-i$ to a vertex in $C_r$, and no shorter such path exists.
Hence $(C_r,\ldots,C_1,C_0)$ is the distance partition with respect to
$C_r$ and so $C_r$ is completely regular.  The result now follows.\qed

Of course, the same argument shows that $C_0$ is completely transitive
if and only if $C_r$ is also.  It will be convenient to denote $C_r$
by $\copp$.  The first author has shown that, in the case where $C$ is
a subset of $J(v,k)$, $C$ and $\copp$ always have the same strength
(see [?????]).  Note that in this case $C_0$ is completely transitive
if and only if the set of complements in $X$ of the $k$-sets in $C_0$
is a completely transitive subset of $J(v,v-k)$.  Thus we may always
assume that $k\le v-k$.

\smallskip\noindent
\eg1
	Let $Y$ be a subset of $X$ and let $C_0$ be the set of $k$-subsets
	$\al$ of $X$ such that $\al\cap Y$ is maximal.  (So if $k\ge|Y|$ then
	$C_0$ consists all $k$-subsets which contain $Y$; otherwise it is all
	the $k$-subsets of $Y$.)  Then $\copp$ consists of the $k$-subsets of
	$X$ whose intersection with $Y$ is minimal.

\smallskip\noindent
In the next four examples, we assume that $\{Y_1,\ldots,Y_b\}$ is a
partition of $X$ with $|Y_i|=a$ for all $i$.

\smallskip\noindent
\eg2 Assume $b=2$ and $k\le a$.  Let $C_0$ be the set of all $k$-sets
contained in $Y_1$ or $Y_2$.  (In this case $\copp$ consists of all
$k$-sets which meet one of the $Y_i$ in $\lfloor k/2\rfloor$ points
and the other in $\lceil k/2\rceil$ points.)

\smallskip\noindent
\eg3 Assume $a=2$ and let $C_0$ be the set of all $k$-sets containing at
most one element from each $Y_i$.  (In this case $\copp$ consists of
the $k$-sets which meet at most one of the $Y_i$ in a single point.)

\smallskip\noindent
\eg4 Assume $a,b\ge3$ and $k=3$.  Let $C_0$ be the set of triples meeting
each $Y_i$ in at most one point.

\smallskip\noindent
\eg5 Assume $k=2$ and let $C_0$ be the set of all pairs meeting each $Y_i$
in at most one point.

\smallskip\noindent
In each of Examples~2--5 the automorphism group of $C_0$ (or $\copp$)
is $\sym{a}\wr\sym{b}$.  In Example~1 the automorphism group is
intransitive.  We will prove that if $C_0$ is intransitive and its
automorphism group is in intransitive, or transitive but imprimitive,
then $C_0$ is one of the above examples.

Next we mention a number of sporadic examples.  The set of lines of a
projective plane of order two is a completely transitive subset of
$J(7,3)$ with covering radius one, while the set of lines of a plane
of order three is completely transitive in $J(13,4)$.  Martin
\cite{mart} shows that no other projective plane is completely
regular.  Delsarte~\cite{del} observed that the Witt design on 24 points is
completely regular, with covering radius two.  From Table~1 in
\cite{krmm} we see that $M_{24}$ has three orbits on sets of size
eight, hence $W_{24}$ is completely transitive.  Martin \cite{mart}
proved that the Witt design on 23 points is completely regular with
covering radius three.  From Table~1 in \cite{krmm} we see that
$M_{23}$ has four orbits on sets of size seven, whence the design is
completely transitive.  By \cite{mart} the Witt design on 22 points
is not even completely regular (Martin \cite{mart,mar3}).

\section{The Intransitive and Imprimitive Cases}

Our first result is a characterisation of the completely transitive
designs for which the automorphism group $G$ is intransitive on $X$.
A group of automorphisms of $J(v,k)$ will be said to be {\sl
completely transitive} if its orbits on $k$-subsets form the
distance partition for some subset of $J(v,k)$.

\begin{lemma}
	Let $C$ be a completely transitive design on a $v$-set
	$X$ with $0<k<v$ admitting a completely transitive group which is
	intransitive on $X$.  Then $C$ is as in Example 1 for some non-empty
	proper subset $Y$ of $X$ and $\aut{C}\cong \sym{Y}\times\sym{\ovy}$.
	Any completely transitive group with $C_0$ as an orbit is transitive
	on $\binom{Y}{j}\times\binom{\ovy}{k-j}$ for each $j$ such that
	$0\le j\le \min\{k, |Y|\}$.
\end{lemma}

\proof
Let $C_0,\ldots,C_r$ be the cells of the distance partition of $C$ and
let $G\le\sym{X}$ be a completely transitive group preserving $C$ and
intransitive on $X$.  If $Y$ is an orbit of $G$ on $X$ then any two
$k$-sets in the same cell $C_i$ must meet $Y$ in the same number of
points.  Suppose $\al\in C$ and $Y$ is an orbit of $\aut{C}$ such that
$|\al\cap Y|=i$, for some positive integer $i$.  Suppose further that
we have elements of $y$, $y'$, $z$ and $z'$ of $X$ such that
\[
y\in\al\cap Y,\quad y'\in Y\diff\al,
	\quad z\in\al\diff Y,\quad z'\in\ovy\diff\al.
\]
Then the $k$-sets $(\al\diff \{ z\} )\cup \{y'\}$ and $(\al\diff\{ y\}
)\cup \{ z'\}$ are both adjacent to $\al$ in $J(v,k)$, but meet $Y$ in
$i+1$ and $i-1$ points respectively.  Therefore neither of these sets
belongs to $C_0$ and they cannot both lie in $C_1$.  Thus we conclude
that one of the following holds:
\begin{enumerate}[(i)]
	\item
	$Y\sbs\al$,
	\item
	$\al\sbs Y$,
	\item
	$\ovy\sbs\al$.
\end{enumerate}
In the second case we may replace the $k$-sets on $C_0$ by their
complements in $X$, thus reducing to case (i).  In the third case, we
may replace $Y$ by an orbit of $\aut{C}$ contained in $\ovy$, with the
same result.

Suppose now that (i) is true and that $|Y|=m$.  We aim to show that
$C_i$ consists of all $k$-subsets of $Y$ which meet $Y$ in $m-i$
points.  Since $k<v$ we may choose a point $z$ not in $\al$ and a
point $y$ in $Y$.  Then $(\al\diff\{ y\} )\cup\{ z\}$ is adjacent to
$\al$ in $J(v,k)$ and meets $Y$ in $m-1$ points.  It follows that this
$k$-set is in $C_1$, and hence that all elements of $C_1$ meet $Y$ in
$m-1$ points.  Thus if $x\in\al\diff Y$ and $z\notin\al$ then
$(\al\diff \{ x\} )\cup\{ z\}$ must lie in $C_0$ and accordingly $C_0$
consists of all the $k$-subsets of $X$ which contain $Y$.  A $k$-set
which contains exactly $m-i$ points from $Y$ meets each element of
$C_0$ in at most $k-i$ points and hence is at distance at least $i$ in
$J(v,k)$ from $C_0$.  From this it follows easily that $C_i$ consists
of all $k$-subsets of $X$ which meet $Y$ in $m-i$ points, as asserted,
and that $\aut{C}$ is transitive on 
$\binom{Y}{j}\times\binom{\ovy}{k-j}$ 
whenever $0 \leq i \leq m$.  Clearly this will also hold for any
subgroup of $\aut{C}$ which is transitive on each of the sets
$C_i$.\qed

Next we shall classify all completely transitive designs admitting a
transitive imprimitive group $G$.

\begin{lemma}
Let $C$ be a subset of $V(J(v,k))$, with $k \leq v/2$
and suppose that some transitive imprimitive subgroup $G$ of $\sym{X}$
acts completely transitively on $C_0$.  Then $C$ is one of
Examples~2-5, or its opposite.
\end{lemma}

\proof
Let $Y=\{Y_1, \ldots,Y_b\}$ be a non-trivial partition of $X$
preserved by $G$, where $|Y_i|=a$ for all $i$ and $v = ab$.  Let
$C_0,\ldots,C_r$ be the distance partition of $C$, where $C_0=C$.  Let
$\al $ be a $k$-subset of $X$ and for each $i = 1,\ldots, b$ set
$\al_i =\al\cap Y_i$ and $e_i =\ |\al_i|$. Let ${\bf e}(\al)$ be the
multiset $\{ e_1,\ldots,e_b\}$.  If $\al$ and $\be$ lie in the same
cell of the distance partition of $C$ then $\bfe\al=\bfe\be$.

We first consider the case $b=2$. Suppose that $\al$ from $C$ meets
both $Y_1$ and $Y_2$ and that
$$
	|\al\cap Y_1|\ge|\al\cap Y_2|.
$$
If $\al\supseteq Y_1$ then, since $k\le v/2=a$, it follows that
$\al=Y_1$ and $C_0=\{Y_1,Y_2\}$.  (This is an uninteresting special
case of Example~2.)  So we may assume that there are points $x\in\al
\cap Y_1$, $y_1\in Y_1\diff\al$ and $y_2\in Y_2\diff\al$.  Define
$k$-subsets of $X$ by
$$
	\be_1:=(\al\diff \{ x\} )\cup\{ y_1\} ,\quad 
	\be_2:=(\al\diff\{ x\} )\cup\{ y_2\} .
$$
Then $\be_1$ and $\be_2$ are both adjacent to $\al$ in $J(v,k)$ but
$\bfe\al=\bfe{\be_1}$, and $\bfe\al\ne\bfe{\be_2}$ unless $|\al\cap
Y_1|=(k+1)/2$. Suppose first that $\al\subset Y_1$. Then $\be_2\in
C_1$ and $\be_1\in C_0$.  It follows that all $k$-subsets of $Y_1$ lie
in $C_0$; so that $C=C_0$ consists of all $k$-subsets of $Y_1$ or
$Y_2$, i.e, it is as in Example~2. Now suppose that $\al\not\subset
Y_1$. Then $k>i:=|\al\cap Y_1|\ge |\al\cap Y_2|=k-i>0$. Thus there is
a point $x_2\in\al\cap Y_2$ and the $k$-subset $\be_3:=(\al\diff\{
x_2\} )\cup\{ y_1\}$ is adjacent to $\al$ in $J(v,k)$. Morever
$\bfe{\be_3}\ne\bfe{\al }$, so $\be_3\in C_1$. So now we have
$\bfe{\be_2}\ne\bfe{\be_3}$ (since $i\ge k/2$) and
$\bfe{\be_1}\ne\bfe{\be_3}$, and hence $\be_1,\be_2\notin C_1$. Since
$\be_1, \be_2$ are both adjacent to $\al$ it follows that
$\be_1,\be_2\in C_0$ whence $\bfe{\be_2}=\bfe{\al}$. Therefore
$i=(k+1)/2$, and $C_0$ consists of all $k$-subsets which contain
$(k+1)/2$ points of one of the $Y_j$ and $(k-1)/2$ points of the
other. In this case $\copp$ consists of all $k$-subsets of $Y_1$ or
$Y_2$ and so $\copp$ is as in Example~2.

Next consider the case $a=2$.  Call a subset of $X$ which meets each
set $Y_i$ in at most one point a \textsl{partial transversal}.  Suppose
$\al\in C$ meets at least two of the sets $Y_i$ in one point, but is
not a partial transversal.  Then we may assume that $\al$ is disjoint
from $Y_1$, meets both $Y_2$ and $Y_3$ in one point, and contains
$Y_4$.  Let $\be_1$ be the $k$-set obtained from $\al$ by deleting one
of the points in $Y_4$ and replacing it with a point in $Y_1$.  Let
$\be_2$ be the $k$-set obtained by deleting the point in $Y_2$ from
$\al$ and replacing it with the point in $Y_3\diff\al$.  Finally let
$\be_3$ be obtained by removing one of the points in $Y_4$ and adding
the missing point from $Y_3$.  Then $\bfe{\be_1}$, $\bfe{\be_2}$ and
$\bfe{\be_3}$ are all distinct.  But if $\al\in C_0$ then all its
neighbours must lie in $C_0\cup C_1$.  Thus we conclude that either
$\al$ is a partial transversal or it must meet at most one of the
$Y_i$ in a single point.

Assume first that $\al$ from $C_0$ is a partial transversal.  Then all
elements of $C_0$ are partial transversals.  The neighbours of $\al$
in $J(v,k)$ are either partial transversals or contain exactly one of
the sets $Y_i$.  From this it follows $C_0$ is set of all partial
transversals of size $k$, and so $C_0$ is as in Example~3.  Thus we
are left with the possibility that $\al$ meets at most one $Y_i$ in a
single point, and thus all elements of $C_0$ contain exactly $\lfloor
k/2\rfloor$ of the sets $Y_i$.  In this case $\copp$ must contain a
partial transversal and so, by what we have just proved, $C_0$ must be
the opposite of a completely transitive subset as in Example~3.

If $k=2$ or $3$ then it is easy to see that $C_0$ is as in Example~4
or Example~5 respectively.  Accordingly we may assume that $a\ge3$,
$b\ge3$ and $k\ge4$.  Let $k_1$ and $a_1$ be respectively the
quotient and remainder when $k$ is divided by $a$.  Then there is a
$k$-subset of $\al$ of $X$ such that $k_1$ of the entries of
$\bfe\al$ are equal to $a$, one is $a_1$ and the rest are all zero.
Note that since $v=ab$ and $k\le v/2$,
\begin{equation}
	\label{eq:k1ka}
	k_1=\lfloor \frac{k}{a}\rfloor \le \lfloor \frac{b}{2}\rfloor
\end{equation}
and consequently there must be at least one entry of $\bfe\al$ equal
to zero.  Assume for the moment that $k_1>0$ and $0<a_1<a-1$.  Then
there is a $k$-subset $\be$ adjacent to $\al$ with $k_1-1$ of the
entries of $\bfe\be$ equal to $a$, one entry equal to $a-1$, another
equal to $a_1+1$, and the remaining entries zero.  There is also a
$k$-subset $\ga$ adjacent to both $\al$ and $\be$ with $k_1-1$ of the
entries of $\bfe\ga$ equal to $a$, one entry equal to each of $a-1$,
$a_1$ and $1$, and the rest zero.  Suppose that $\al\in C_i$ for some
$i$.  Since $\bfe\al$, $\bfe\be$ and $\bfe\ga$ are all different, one
of $\be$ and $\ga$ must lie in $C_{i-1}$ and the other in $C_{i+1}$.
As $\be$ and $\ga$ are adjacent, this is impossible.

We are left with the possibilities that $k_1=0$ or that $a_1=0$ or
$a-1$.  An obvious analogue of the argument of the last paragraph will
work if we have a $k$-subset $\al$ such that $\bfe\al$ has at least
one entry equal to each of $a$, $x$ and $0$, where $a-2\ge x>0$.
Suppose first that $a_1=0$.  Then $0<k_1\le b-2$, by \eqref{eq:k1ka}, and it
follows that there is a $k$-subset $\al$ such that $\bfe\al$ has
$k_1-1$ entries equal to $a$, one entry equal to each of $a-1$ and
$1$, and the rest zero.  If $k_1\ge 2$ then by \eqref{eq:k1ka}, $b\ge k_1+2$,
and we may use $\al$ in the above argument (with $x=1$).  On the other
hand if $k_1=1$ so that $k=a$, then there are mutually adjacent
$k$-subsets $\al ,\be ,\ga$ with $\bfe{\al}=\{ a-1,1,0^{b-2}\}$,
$\bfe{\be}=\{ a-2,2,0^{b-2}\}$, and $\bfe{\ga}=\{ a-2,1,1,
0^{b-3}\}$. Next suppose that $a_1=a-1$. In this case if $k_1>0$ then
there is a $k$-subset $\al$ with $k_1$ of the entries of $\bfe\al$
equal to $a$, one equal to each of $a-2$ and $1$, and the rest zero.
We may use $\al$ in the above argument (with $x=a-2$) provided that
$b>k_1+2$, while if $b\le k_1+2$ then by \eqref{eq:k1ka}, $b$ is 3 or 4 and
$k=(k_1+1)a-1\ge (b-1)a-1>v/2$ which is not allowed. This leaves the
case $k_1=0$, $k=a_1$. In this case there are mutually adjacent
$k$-subsets $\al ,\be ,\ga$ with $\bfe{\al}=\{ k-1,1,0^{b-2}\}$,
$\bfe{\be}=\{ k-2,2,0^{b-2}\}$, and $\bfe{\ga}=\{ k-2,1,1,
0^{b-3}\}$.\qed

The problem of classifying completely transitive designs admitting
primitive groups is open, but by the results of this section is
reduced to the following problem.

\begin{problem}
Classify all completely transitive designs on a $v$-set $X$
admitting a group $G \leq \sym{X}$ which is primitive on $X$.
\end{problem}

\section{Designs with Minimum Distance at Least Three}

As in the previous section $C$ is a completely transitive design in
$J(v,k)$.  If $|C| = 1$ then it is a case of Example 1.1 with $\ell =
k$, so assume that $|C|\ge2$.  Let $\de$ denote the minimum distance
between two $k$-subsets in $C$.  The design of Examples~2.1 has
$\de=1$, hence Proposition~3.1 yields that the automorphism group of a
comletely transitive design with $\de\ge2$ must be transitive on the
underlying set $X$.  Most of the designs in Examples~2.2--5 have
$\de=1$.  We describe the exceptions explicitly.  The first arises
when $v=2k$, when we may take $C$ to consist of two disjoint
$k$-subsets of $X$.  (This is a special case of Example~2.2, and has
$\de=k$.)  If $v$ and $k$ are even and $C$ consists of the $k$-subsets
formed by the union of any $k/2$ of the cells of a fixed partition of
$X$ into pairs, then $C$ is completely transitive with $\de=2$.  If
$v$ is divisible by three then the $v/3$ triples in a fixed partition
of $X$ into 3-sets is a completely transitive design with $k=\de=3$.

Meyerowitz \cite{mey} shows that the completely regular subsets of
$J(v,k)$ with strength zero are precisely the designs of Example~2.1
(and their opposites).  We note that Martin \cite{mart,mar2} has
shown that a completely regular subset of $J(v,k)$ with strength one
and $\de\ge2$ must be one of the designs just described above.
Consequently any completely regular subset of $J(v,k)$ with $\de\ge3$
must be have strength at least two, and is therefore a 2-design in the
usual sense of the word. From now on we concentrate on the case
$\de\ge3$.

\begin{theorem}
If $C$ is a completely transitive design with $|C|\ge2$ and
$\de\ge3$ and $\al\in C$ then $\aut{C}_\al$ is transitive on the
Cartesian product $\al\times(X\bs\al)$.  Further, either:
\begin{enumerate}[(a)]
	\item
	$v=2k\ge6$ and $C$ consists of two disjoint $k$-subsets,
	\item
	$v=3b$ and $C$ consists of $b$ pairwise disjoint triples, or
	\item
	$\aut{C}$ is 2-transitive on $X$.
\end{enumerate}
\end{theorem}

\proof
Let $C_0,\ldots,C_r$ be the distance partition of $C$, where $C_0=C$,
and assume $\al\in C_0$.  Assume $x,x' \in \al$ and $y, y' \in X\bs
\al$.  Then $\beta = (\al\bs\{x\}) \cup\{y\}$ and $\beta'
=(\al\bs\{x'\})\cup \{y'\}$ are both adjacent to $\al$ and, as
$\de\ge3$, both $\beta$ and $\beta'$ lie in $C_1$.  Hence $\beta^g
=\beta'$ for some $g \in G$.  If $\al^g\ne\al$ then we $\al^g \in
C^g_0 = C_0$ and so $\al$ and $\al^g$ are at distance at most two in
$J(v,k)$, which is a contradiction.  Hence $\al^g = \al$, so $g\in
G_\al $ and $g$ maps $x$ to $x'$ and $y$ to $y'$.  Further, if we take
$x = x'$ above then we see that $\aut{C}_{\al,x}$ is transitive on
$X\bs\al$, whence $\aut{C}_\al$ is transitive on $\al \times
(X\bs\al)$.

By Lemma 3.1 either $C$ is as in the statement or $\aut{C}$ is
primitive on $X$.  Assume the latter holds, and let $B(x)$ be the
intersection of all $k$-subsets of $C_0$ which contain $x$.  It is
easy to see that $B(x)$ is a block of imprimitivity for $\aut{C}$ in
$X$ containing $x$, and as $\aut{C}$ is primitive we have
$B(x)=\{x\}$.  Let $y, y'$ be distinct points of $X\diff\{x\}$. Then
there are $\al, \al' \in C_0$ containing $x$, such that $y \notin \al$
and $y' \notin \al'$.  Now $X\bs \al$ and $X\bs \al'$ have a point in
common, $z$ say, since $k\le v/2$ and $x \in \al\cap\al'$.  Since
$\aut{C}_{\al,x}$ and $\aut{C}_{\al' ,x}$ are transitive on $X\bs \al$
and $X\bs \al'$ respectively there are elements $g$ in $\aut{C}_{\al ,
x}$ and $g'$ in $\aut{C}_{\al', x}$ such that $y^g = z$ and $z^{g'} =
y'$. It follows that $\aut{C}_x$ is transitive on $X\diff\{x\}$, that
is $\aut{C}$ is 2-transitive on $X$.\qed

Next we deal with the case where $k$ is small.

\begin{theorem}
Let $C$ be a completely transitive design on a $v$-set $X$
with $k\le v/2$ such that $|C|\ge2$, $\de\ge3$ and $G$ is
2-transitive.  Then either
\begin{enumerate}[(a)]
	\item
	$v=13$, $k = 4$, $r = 2$, and $C_0$ is the set of lines of 
	the Pappian projective plane $PG_2(3)$, or
	\item
	$k\ge5$ and $r\le k-2$.
\end{enumerate}
\end{theorem}

\proof
Since $\aut{C}$ is 2-transitive on $X$, for any $k$-subset $\be$ and
any pair of points $\{y,y'\} \sbs \be$, there is some $\al \in C_0$
containing $\{y,y'\}$.  Hence $\al$ and $\be$ are distance at most
$k-2$ in $J(v,k)$ whence $\be\in C_i$ for some $i\le k-2$. Hence $r\le
k-2$. Since $\de\ge3$, we have $k\ge \de\ge3$, and $r\ge1$.  If $k =
3$ then, for $\beta = \{x, y, z\} \in C_1$, there are points $s$ and
$t$ in $X\diff\be$ such that $\al = \{x, y, s\}$ and $\al' = \{x, z,
t\}$ belong to $C_0$ (since $G$ is 2-transitive on $X$); then $\al$
and $\al'$ are at distance two in $J(v,k)$ contradicting $\de\geq
3$. Hence $k\ge4$.

Let $k=4$. If $\al = \{x, y, z, w\}$ and $\al' = \{x, y, z', w'\}$ are
distinct 4-subsets in $C_0$ both containing $x$ and $y$ then $\al$ and
$\al'$ are at distance at most two in $J(v,k)$, contradicting $\de\geq
3$. Hence each pair from $X$ lies in a unique 4-subset in $C_0$.  So
$C_0$ is the set of blocks of a 2-transitive $2-(v,4,1)$ design.
There are $\binom{v-2}{2}$ subsets of size four containing $x$ and
$y$, one of which is in $C _0$ and $2(v - 2)$ of which are in $C_1$.
Since $v \ge2k = 8$ it follows that $\binom{v-2}{2}>2(v-2)+1$, so $C_2
\not= \emptyset$.  As any 4-subset contains at least two points of
some 4-subset in $C_0$ it follows that $r=2$, and $\aut{C}$ is
transitive on independent 4-subsets, that is 4-subsets in which no
three points lie in a block of $C_0$.

By \cite{kan}, the 2-transitive $2-(v,4,1)$ designs are the projective
spaces $PG_d(3)$ and affine spaces $AG_d(4)$ with $d\ge 2$, and the
Hermitian and Ree unitals on 28 points. Since $PGL_3(3)$ is transitive
on independent 4-subsets of points, and also on triples $(l,x,y)$,
where the point $x$ lies on the line $l$ and the point $y$ does not,
it follows that taking $C_0$ to be the set of lines of $PG_2(3)$ gives
an example. If $d\ge 3$ then independent 4-subsets might span a plane
or a three-dimensional projective space, so there are no further
examples from $PG_d(3)$.  Similarly for $AG_d(4)$, if $d\ge 3$ then
independent 4-subsets might span a plane or a three-dimensional affine
space, so $d=2$.  But for $AG_2(4)$, $|C_2|= 2^3.3.5.7$, and as 7 does
not divide $|A\Ga L_2(4)$, $\aut{C}$ cannot be transitive on
$C_2$. For the Hermitian and Ree unitals, $|C_0|=63$,
$|C_1|=|C_0|.4.24$, and hence $|C_2|=\binom{28}{4}-97.|C_0|$, which is
greater than $|\aut{C}|$ in either case, so $\aut{C}$ cannot be
transitive on $C_2$.\qed

\begin{corollary}
	\label{cor:rk2}
	If $r\le k-2$ then $|\aut{C}|\ge\frac{1}{k-1}\binom{v}{k}$.
	Unless $k=6$ and $v=12$, this is at least $\frac{1}{4}\binom{v}{5}$.
\end{corollary}

\proof
The group $\aut{C}$ has at most $k-1$ orbits on $k$-subsets from $X$,
so some $\aut{C}$-orbit on $k$-subsets has length at least 
$\frac{1}{k-1}\binom{v}{k}$. A simple arithmetic computation shows that this is
at least $\frac{1}{4}\binom{v}{5}$ unless $k = 6, v = 12$.\qed

In fact if $k < v/2$ then we have
$$
|\aut{C}|\ \geq \frac{1}{k-1}\binom{v}{k} \ge\frac{1}{k-2}\binom{v}{k-1} 
	\ge\cdots\ge\frac{1}{4}\binom{v}{5}.
$$
Applying Corollary~4.3 to some of the known 2-transitive groups is
surprisingly effective.  We prove the following Lemma as a sample
result. 

\begin{lemma}  
	If $r\le k-2$ then the socle $T$ of $G$ is not one of
	$L_2(q)$, $Sz(q)$, $U_3(q)$, $Ree(q)$, for any $q$. 
\end{lemma}
	
\proof
We deal with the four families separately.
\medbreak
\noindent {\bf Suzuki groups:}\quad 
$Sz(q) \le G \le Aut\ Sz(q)$, $q=2^{2s+1}$ odd, $s\ge1$.  Here
$\frac{1}{4} \binom{v}{5} \leq\ | Aut\ Sz(q)|$ becomes
$$
\frac{1}{4} \binom{q^2+1}{5} \leq (q^2 + 1)q^2 (q - 1)a
$$
which implies $q^5 <(q + 1)(q^2 - 2)(q^2 - 3) \leq 480a$, which is not true.
\medbreak
\noindent {\bf Unitary groups:}\quad $U_3(q)\le G \le P\Ga U_3(q)$,
with $q = p^a$ for some prime $p$, and $v = q^3 + 1$. Here the inequality is
$$
\frac{1}{4}\binom{q^3+1}{5} \leq (q^3 + 1)\cdot q^3\cdot (q^2 - 1)\cdot2a
$$
which implies
$(q - 1)^7 < \frac{(q^3 - 1)(q^3 - 2)(q^3 - 3)}{(q^2-1)} \leq 960a$. The
inequality $(q - 1)^7 < 960$ implies $a = 1, q = 2$ or $3$. The exact
inequality above holds only
for $q = 2$, but then $k \ge5 >v/2$.

\medbreak
\noindent {\bf Ree groups:}\quad 
$Ree(q)\le G\le Aut(Ree(q))$, $q=3^{2s+1}$ and $v = q^3 + 1$. The
arithmetic is similar to and easier than that for the unitary groups
and shows that the inequality is never satisfied.

\medbreak
\noindent{\bf 1-dimensional linear groups:}\quad
$L_2(q) \le G \le P\Ga L_2(q)$, $q = p^a$ for some prime $p$ and $v =
q + 1 \geq 2k \geq 10$. If $q = 9$ then $k = 5$ and by Theorem~4.1,
$25 = k(v - k)$ divides $|G|$ which is not the case. Similarly if $q
=11$ then $k$ is five or six and $k(v - k)$ does not divide $|
G_\al|$. Hence $v = q + 1 \ge14$ and by Corollary~4.3,
$$
\frac{1}{4}\binom{q + 1}{5} \leq (q + 1) q(q - 1) a
$$
which implies $(q - 2)(q - 3) \leq 480a$, whence $a \leq 5$. We find
from this that $q\in \{32, 27, 25, 23, 19, 17, 16, 13\}$.  If $k=5$
then using the fact that $5(v - 5)$ divides $| G|$ we obtain $q = 16$
as the only possibility.  Then the 5-subsets $\al$ in $C_0$ are orbits
of subgroups $Z_5$ of $L_2(16)$, and all such orbits form a single
$L_2(16)$-orbit on 5-subsets.  Thus $C_0$ is uniquely determined by
$L_2(16)$, and it follows that $C_0$ is the set of blocks (circles) of
the Miquelian inversive plane of order four, which is a $3-(17,5,1)$
design. Since any triple of points lies on a unique block in $C_0$, it
follows that any two 5-subsets of points are at distance at most 2,
and hence that $r=2$.  However, we have in this case that $|C_0|=68,
|C_1|=|C_0|.60$, and hence that
$|C_2| =\binom{17}{5}-|C_0|-|C_1|=2^2.11^2.17$, which does not divide
$|G|$, and so $G$ is not transitive on $C_2$. Hence $k\ge6$.  Thus we
have $\frac{1}{5}\binom{q+1}{6}\le\ | G|$ whence $(q - 2)(q -3)(q -
4)\le 3600a$, which implies $q \in \{17, 16, 13\}$. Then as $k(v - k)$
divides $| G|$ we must have $q = 17, k = 6$; but as three does not
divide $| G_x|$ it is not possible for $G_{\al, x}$ to be transitive
on $X\diff\al$ where $x \in\al\in C_0$.

\medskip
Next we partially analyse the situation for 2-transitive groups which do not 
fit into any infinite family of 2-transitive groups, namely the Mathieu
groups, $L_2(11)$ of degree 11, and the Higman-Sims and Conway groups $HS$ and $Co_2$. 


\medbreak\noindent{\bf Mathieu groups}\quad
$M_v$ where $v\in\{11, 12, 22, 23, 24\}$, or $L_2(11)$ with $v=11$, or
$M_{11}$ with $v = 12:$\quad Since $k(v - k)$ divides $| G_\al|$ we
have the following possibilities:
\begin{enumerate}[(i)]
	\item
	$v = 11$, $k = 5$, $G = M_{11}$ or $L_2(11)$,
	\item
	$v = 12$, $k = 6$, $G = M_{12}$ or $M_{11}$,
	\item
	$v = 22$, $6\le k\le10$, $k\ne9$, $G = M_{22}$ or $\aut{M_{22}}$
	\item
	$v = 23$, $5\le k\le11$, $k\ne6, 10$, $G = M_{23}$
	\item
	$v = 24$, $6\le k\le12$, $k\ne7, 11$, $G = M_{24}.$
\end{enumerate}
We consider the cases separately. In case (i), if $G=M_{11}$, then $G$
has two orbits on 5-subsets of points and we obtain two completely
transitive designs, namely the Witt design on 11 points and its
opposite, both preserved by $M_{11}$.  However in this case, $\de=2$
:-(. If $G=L_2(11)$ then $G$ has 4 orbits on 5-subsets of points.
Let $C_0$ be the set of blocks of the $2-(11,5,2)$ design
preserved by $G$.  Then this is one of the orbits.  Since each pair of
blocks in $C_0$ intersect in 2 points, $\de=3$.  If $\al\in C_0$ then,
since $G_\al$ is transitive on $\al\times(X\diff\al )$, it follows
that $G$ is transitive on the 330 subsets in $C_1$.  Next suppose
$\al'\in C_0\diff\{\al\}$.  Then the only 5-subset containing
$\al\cap\al'$ and at distance three from $\al$ is $(\al\cap\al' )\cup
(X\diff\al\cup\al')$, and this is fixed setwise by $G_{\al\cap\al' }$;
since $G_{\al\cap\al' }$ is maximal in $G$ it follows that these
5-subsets form a $G$-orbit $C_3$ of size 55.  Finally, a 5-subset
$\beta$ contained in $X\diff\al$ has stabiliser $D_{10}$, and we
deduce that all 5-subsets contained in the complement of a block of
$C_0$ form a $G$-orbit of length 66 which must be $C_2$.  The order of
the stabiliser of a 5-subset in $C_1$, $C_2$ and $C_3$ is $2$, $12$
and $10$ respectively so, by Theorem~4.1, none of these orbits can be
completely transitive with $\de\ge3$.  The covering radius of $C_0$ is
two, not three, and so $C_0$ is also not completely transitive.
(Presumably none of these orbits is even completely regular, this is
certainly the case for $C_0$ [??].)

Now consider case (ii). If $G=M_{12}$, then there are two orbits on
6-subsets, we obtain a completely transitive design, but again in this
case we have $\de=2$.  Thus $G=M_{11}$.  By \cite{atlas}, the inner
product of the permutation characters for the actions of $M_{12}$ on
the cosets of this subgroup $G$ and on the blocks of the $5-(12,6,1)$
Steiner system is two, whence $G$ has two orbits, $C_0$ and $C_3$ say,
on the set $B$ of 132 blocks of the $5-(12,6,1)$ Steiner system. These
orbits have lengths 22 and 110, the block stabilisers being $A_6$ and
$3^2:[2^4]$ respectively.  Since both of these subgroups are
transitive on $\al\times(X\diff\al )$, where $\al\in B$ is fixed by
the subgroup, it follows that $G$ has two orbits on the 6-subsets not
in $B$, namely $C_1$ of length 132, and $C_2$ of length 660.
Moreover, $C_0$ is a completely transitive design, as is its opposite.
Now $C_0$ is the set of blocks of the $3-(12,6,2)$ design preserved by
$G$; since each 3-subset is contained in two blocks of $C_0$ it
follows that $\de=3$.  A computer calculation shows that the dual
degree of this design is two, whence its covering radius is at most
two.  (In fact $C_0$ is completely regular, but $C_2$ is the union of
two orbits of $M_{12}$.)  [Thanks to Bill Martin for the compuations.]

From Table~1 in \cite{krmm} we find that $M_{22}$ has at least $k$
orbits on $k$-sets when $k\ge6$.  Diagram~3 of \cite{krmm} gives
detailed information about the action of $M_{22}$ on subsets of $X$,
and this provides information about $\aut{M_{22}}$ as well.  For any
orbit of $\aut{M_{22}}$ is either an orbit of $M_{22}$, or the union
of two $M_{22}$ orbits of the same length.  From Diagram~3 we see that
if $k\ge6$, there are three $M_{22}$-orbits of different lengths which
are pairwise adjacent in $J(22,k)$.  Thus no completely transitive
designs arise in connection with $M_{22}$ or $\aut{M_{22}}$.  From
Diagram~2 of \cite{krmm} we see that if $8\le k\le12$ then there are
three orbits of $M_{23}$ on $k$-sets, pairwise at distance one in
$J(23,k)$.  So $k\le7$ in this case.  From Diagram~1 of the same
source, we see that $M_{24}$ has three pairwise adjacent orbits on
12-subsets, whence we have $k\le11$.

In cases (i) and (ii) we get an example: $C_0$ is the set of blocks of
the Steiner system.  Similarly in case (iv) with $k = 7$ and case (v)
with $k = 8$ the blocks of the Steiner system give examples.  The Witt
design on 23 points is completely transitive.  [*** We need to
consider the orbits of $M_{23}$ on 5- and 6-subsets; these may
completely transitive, but have $\de<3$?  For $M_{24}$ we have to
consider the set of all $k$-subsets which contain a block of the Witt
design on 24 points.  By \cite{krmm}, $M_{24}$ has three orbits on
$k$-subsets when $8\le k\le11$.  These have $\de<3$ as well.***]

\medbreak\noindent
{\bf Higman Sims and the Conway group $Co_3$:}\quad For $G = HS, v =
176$ and $k(176 - k)$ divides $| G| = 2^9\cdot3^2\cdot5^3\cdot
7\cdot11$. From the inequality $\frac{1}{k - 1}\binom{176}{k} \leq | G|
$ it follows that $k \leq 19$, and then the divisibility condition
implies that $k$ is $8$, $11$ or $16$.
\bigskip
These need to be analysed.  [*** Mmm ***]

\bigbreak
For $G=Co_3$, it follows from Corollary~3.4 that $k\le 6$, and from the
fact that $k(276-k)$ divides $|G|$, we have $k=6$.

\medbreak
The remaining infinite families of 2-transitive groups are the following:
\begin{enumerate}
 \item[(i)] The projective groups:  $L_d(q)\le G \le P\Ga L_d(q)$ 
where $d\ge3$ and $v =(q^d-1)/(q-1)$, or
$G=A_7<L_4(2)$ with $v=15$.

\item[(ii)] The affine groups: $G=N.G_0 \le A\Ga L_d(q)$ where $d\ge1$
and $v=q^d$. Here $N$ is the group of translations acting 
regularly on the vector space $X=F+q^d$,  and $G_0\leq \Ga L_d(q)$
ats transitively on the non-zero vectors. 

\item[(iii)] The Symplectic groups: $G= Sp_{2m}(2)$ acting 2-transitively on
the set of  $v = 2^{2m-1} +\varepsilon 2^{m-1}$ nondegenerate quadratic forms
of type $\varepsilon=\pm$, which polarise to the symplectic form preserved by $G$.  
\end{enumerate}

We do not treat the symplectic groups at all. We make a partial 
analysis of the other two cases.

\medbreak\noindent
{\bf Projective groups:}\quad
$L_d(q)\le G \le P\Ga L_d(q)$ where $d\ge3$ and $v =(q^d-1)/(q-1)$, or
$G=A_7<L_4(2)$ with $v=15$.

There are two constraints on $k$ which must be satisfied if $G$ is
completely transitive.  First, since $G\le P\Ga L_d(q)$, we have
\begin{equation}
	\label{eq:pgaml}
	|P\Ga L_d(q)|\ge \frac{1}{k-1}\binom{v}{k}.
\end{equation}
This provides an upper bound on $k$.  To establish a lower bound, we
first show that $\al$ cannot be an  $i$-dimensional subspace, 
for $1\leq i\leq d-2$.

The size of the symmetric difference of any two $i$-dimensional
subspaces is at least $2q^i$ and thus the minimum distance in $J(v,k)$
between two $i$-dimensional subspaces is $q^i$.  By Theorem~4.2(b) we
may assume that $k\ge5$, so $q^{i+1}-1\ge5(q-1)$ and therefore
$q^i\ge4$.  It follows that if we delete any two points $p$ and $p'$
from $\al$ and replace them by two points $q$ and $q'$ not in $\al$,
the resulting $k$-set $\be$ is not a subspace.  Hence $\be\notin
C_0\cup C_1$, whence it must lie in $C_2$.  Now if $q'$ is not in the
span of $\be\cup q$ then $\be$ spans an $(i+2)$-dimensional subspace,
otherwise it spans an $(i+1)$-dimensional subspace. If $i\leq d-3$ then
both cases may occur, it follows that $G$ has at least two orbits on $C_2$, and is
not completely transitive. 

Suppose then that $i=d-2$, that is, $\alpha$ is 
a hyperplane. Let
$p$ and $p'$ be distinct points from $\al$.  Let $q$ and $r$ be two
points not in $\al$ such that the unique line through $q$ and $r$
meets $\al$ in a point distinct from $p$ and $p'$, and let $q'$ be a
point on the line through $p$ and $q$ which is distinct from $p$ and
$q$.  (So $q'\notin\al$.)  Suppose
$$
\be:=\al\diff\{p,p'\}\cup\{q,r\}\quad \be':=\al\diff\{p,p'\}\cup\{q,q'\}.
$$
Arguing as in the previous paragraph, we see that both $\be$ and
$\be'$ lie in $C_2$.  Suppose that there is an element $g$ of $G$ such
that $\be^g=\be'$.

If $\al^g$ is not equal to $\al$ then we have $|\al\cap\al'|\ge k-4$,
but $\al\cap\al'$ is a subspace and therefore has size at most
$k-q^{d-2}$.  So $q^{d-2}\le4$, that is either $G = L_4(2)$ and $k =
7$ or, since $k\ge5$, we have $G\ge L_3(4)$ and $k = 5$.  In both
cases $|\al\cap\al^g|=k-4$ and $\de=4$, and so $\al^g$ contains $\{q,
q'\}$ and $\al\cap\al^g\sbs\al\diff\{p, p'\}$.  In the case where $G
\ge L_3(4)$ this means that the line containing $\{q,q'\}$ is $\al^g$
and it meets $\al$ in a point different from $x$, which is a
contradiction.  Similarly in the case $G=L_4(2)$ the line $\{p, q,
q'\}$ must lie in the subspace $\al'$ which is not the case.  Hence
$\al^g=\al$, and so $(\al\diff\{x, x'\})^g \sbs\al\cap\be'$ and
consequently $(\al\diff\{p, p'\})^g = \al\diff\{p, p'\}$.  Therefore
$\{p,q,q'\}^{g^{-1}} =\{p^{g^{-1}}, q, r\}$ is a collinear triple,
which is a contradiction since $p^{g^{-1}} \in \{p, p'\}$.  Thus we
have shown that $\al$ spans the whole space $X$.

If $\al$ is not a subspace it follows that there is a line which meets
$\al$ in at least two points, but is not contained in it.  Suppose $p$
is a point in $\al$ and $\ell$ is a line on $p$ which meets $\al$ in
exactly $x$ points, where $2\le x\le q$.  By Theorem~4.1 we know that
$G_{\al,p}$ is transitive on those lines through $p$ which are not
contained in $\al$.  Hence every line on $p$ meets $\al$ in at least
$x$ points.  It follows that the number of lines on $p$ is a lower
bound on $k-1$, and thus we have
\begin{equation}
	\label{eq:k1v1}
	k-1\ge \frac{v-1}{q}.
\end{equation}

We are now going to apply this to the inequality in \eqref{eq:pgaml}.  If $p$
is prime and $q=p^a$ then
\begin{equation}\label{order}
|P\Ga L(d,q)| =\frac{a}{q-1}(q^d-1)\cdots(q^d-q^{d-1})< q^{d^2}.
 \end{equation}
We need to compare this with the ratio $\binom{v}{k}/(k-1)$.  A
routine calculation shows that when $v\ge2k$, this is an increasing
function of $k$.  Since $v/k<(v-i)/(k-i)$ when $0<i<k$ we also have
$$
\binom{v}{k}> \left(\frac{v}{k}\right)^k.
$$
From \eqref{eq:k1v1} we see that $k\ge v/q$, whence $v/k\ge q$ and $k\ge
q^{d-2}+q^{d-3}$.  (Note that $d\ge3$.)  Accordingly \eqref{order}
implies that if $G$ is completely transitive then
$$
q^{d^2+d-1}\ge (k-1)q^{d^2}> q^{q^{d-2}+q^{d-3}},
$$
implying in turn that 
\begin{equation}
	\label{eq:d2d1}
	d^2+d-1>q^{d-2}+q^{d-3}.
\end{equation}
This yields the following possibilities:
$$
d=3,\ q\le9;\qquad d=4,\ q\le3;\qquad 5\le d\le 7,\ q=2.
$$
Using \eqref{eq:k1v1} with $k=\lfloor (v-1)/q \rfloor+1$ in place of
\eqref{eq:d2d1} the above list reduces to:
$$
d=3,\ q\le5;\qquad d=4,\ q=2.
$$
[*** So here are six more cases; I know how to eliminate the $q=2$ ones ***]

\medbreak\noindent
{\bf Affine groups:}\quad
$G=N.G_0 \le A\Ga L_d(q)$ where $d\ge2$, $v=q^d$, $N$
is the group of translations and $G_0\le\Ga L_d(q)$ is 
transitive on non-zero vectors. 

We aim to show that if $G$ is completely transitive then $k\ge
q^{d-1}$.  Suppose $\al\in C_0$ and $\al$ is an affine subspace with
dimension $i$, where $i\le d-2$.  The size of the symmetric difference
of any two $i$-dimensional subspaces is at least $2(q-1)q^{i-1}$.  By
Theorem~4.2(b) we have that $q^i=k\ge5$ and so any two $i$-dimensional
subspaces are at distance at least four in $J(v,k)$.  As in the
projective case, it follows that $C_2$ splits into at least two orbits
under the action of $G$.  Thus if $\al$ is a subspace then it has
dimension $d-1$ and $k=q^{d-1}$.

Suppose then that $\al$ is not a subspace and let $p$ be a point in
$\al$.  Since $G_\al$ is transitive on $\al$, there must be a line on
$p$ which contains both a point in $\al\diff p$ and a point not in
$\al$.  As $G_\al$ is transitive on the set of lines through $p$ which
contain a point not in $\al$, it follows that every line through $p$
contains a point of $\al\diff p$.  Therefore
$$
|\al\diff p|=k-1\ge \frac{v-1}{q-1} =\frac{q^d-1}{q-1}> q^{d-1}.
$$
Thus we have shown that $k\geq q^{d-1}$.

If $p$ is prime and $q=p^a$ then
$$
|A\Ga L(d,q)| = a q^d |GL(d,q)| \le q^{d^2+d+1}
$$
while
$$
\binom{v}{k}> \binom{v}{q^{d-1}} > q^{q^{d-1}}.
$$
Consequently, if $G$ is completely transitive, we must have
$$
q^{(q+1)^2}\ge (k-1)q^{d^2+d+1} \ge \binom{v}{k} > q^{q^{d-1}}
$$
and thus
$$
(d+1)^2 > q^{d-1}.
$$
This leaves the possibilities
$$
d=2,\ q\le 8;\qquad d=3,\ q\le3;\qquad 4\le d\le6.
$$
[ *** Presumably some of these can be eliminated by using the exact
values for $|A\Ga L(d,q)|$ and $\binom{v}{k}/(k-1)$; I have not done
these computations yet. ***]

\end{document}